\documentclass[letterpaper, 10 pt, conference]{ieeeconf}  

\IEEEoverridecommandlockouts                              

\usepackage{graphicx,amsmath}
\usepackage{amssymb}
\usepackage{amsfonts}
\usepackage{calrsfs}
\usepackage{mathrsfs}
\usepackage{caption}
\usepackage{subcaption}
\usepackage{multirow}
\usepackage{xcolor}
\usepackage{graphics}
\usepackage{float}
\usepackage{graphicx}
\usepackage{calrsfs}
\usepackage{mathtools}
\usepackage{amsmath}
\usepackage{psfrag}
\usepackage{blindtext}
 \newtheorem{teorema}{{\em Theorem}}{ }
 \newtheorem{lema}{{\em Lemma}}{ }
 { }

\makeatletter
\renewcommand*\env@matrix[1][*\c@MaxMatrixCols c]{%
	\hskip -\arraycolsep
	\let\@ifnextchar\new@ifnextchar
	\array{#1}}
	\makeatother
\mathtoolsset{showonlyrefs}
 \usepackage{graphicx}
\usepackage{amssymb}
\usepackage{epstopdf}
\usepackage{xcolor,colortbl}
\usepackage{hyperref}
\usepackage{pbox}
\usepackage{epsfig,psfrag}
\usepackage{color}
\usepackage{amsmath}
\usepackage{mathrsfs}  
\usepackage{amsbsy}
\usepackage{lipsum} 
\usepackage{pgf}
\usepackage{tikz}
\usepackage{ragged2e}
\usetikzlibrary{arrows, automata, shapes.multipart, chains, positioning, fit, graphs, spy, shapes, backgrounds, patterns,decorations.pathmorphing}
\usepackage{verbatim}
\tikzstyle{block} = [draw, fill=blue!20, rectangle, 
minimum height=1 cm, minimum width=2cm, rounded corners=0.1cm]
\tikzstyle{block2} = [draw, fill=black!15, rectangle, 
minimum height=0.75 cm, minimum width=1.5cm, rounded corners=0.1cm]
\tikzstyle{faultyblock} = [draw, fill=red!20, rectangle, 
minimum height=1 cm, minimum width=2cm, rounded corners=0.1cm]
\tikzstyle{reconblock} = [draw, fill=green!20, rectangle, 
minimum height=1 cm, minimum width=2cm, rounded corners=0.1cm]
\tikzstyle{sum} = [draw, fill=black!15, circle, minimum width=0.5cm]
\tikzstyle{input} = [coordinate]
\tikzstyle{output} = [coordinate]
\tikzstyle{pinstyle} = [pin edge={to-,thin,black}]
\usepackage{float}
\usepackage{isomath}
\usepackage{subcaption}

\title{\LARGE \bf
Extremum Seeking Control for Multivariable Maps\\ under Actuator Saturation
}

\author{Enzo Ferreira Tomaz Silva$^{1}$, Pedro Henrique Silva Coutinho$^2$, Tiago Roux Oliveira$^2$, Miroslav Krsti\'{c}$^3$
\thanks{*This work was supported by the Brazilian agencies CNPq (Grant numbers: 407885/2023-4 and 309008/2022-0), CAPES, and FAPERJ.}
\thanks{$^{1}$ Enzo Silva is with the Graduate Program in Electronics Engineering, State University of Rio de Janeiro, Brazil.
			        {\tt\small enzotomazsilva@gmail.com} }
\thanks{$^{2}$ Pedro Coutinho and Tiago Roux Oliveira are with the Department of Electronics and Telecommunication
Engineering, State University of Rio de Janeiro, Brazil.
			        {\tt\small phcoutinho@eng.uerj.br, tiagoroux@uerj.br }}
\thanks{$^{3}$ M. Krsti\'{c} is with the  Department of Mechanical and Aerospace Engineering, University of California - San Diego, USA (e-mail: krstic@ucsd.edu).}
}

\begin{document}

\maketitle
\thispagestyle{empty}

\pagestyle{empty}

\begin{abstract}
This paper deals with the gradient-based extremum seeking control for multivariable maps under actuator saturation. By exploiting a polytopic embedding of the unknown Hessian, we derive a LMI-based synthesis condition to ensure that the origin of the average closed-loop error system is exponentially stable. Then, the convergence of the extremum seeking control system under actuator saturation to the unknown optimal point is proved by employing Lyapunov stability and averaging theories. Numerical simulations illustrate the efficacy of the proposed approach.
\end{abstract}

\section{Introduction}

With the technological advances of recent decades, real-time optimization methods have been employed for solving several practical problems. Within this context, several strategies of extremum seeking control have been developed and applied in various scenarios~\cite{ESsurvey}.
The extremum seeking control is an adaptive, real-time, and model-free strategy. The purpose of this technique is to find an optimal point such that a given desired function (with unknown parameters) is maximized or minimized, i.e., its extremum point is reached~\cite{LivroTiagoES}.

One of the well-known extremum-seeking approaches is the gradient algorithm. In this case, the optimization process works by applying sinusoidal disturbances in the control scheme, and the algorithm evaluates the gradient direction to adjust the control signal that will forward the system to the optimal point desired.

Several efforts have been made to extend the extremum seeking control to different classes of maps and control problems, such as delay systems~\cite{ArtDelay2}, cascade maps with partial differential equations~\cite{ArtPDE, galvao2022extremum,silva2023extremum}, cooperative games with Nash equilibrium~\cite{rodrigues2024sliding, rodrigues2024nash}, and event-triggered control~\cite{rodrigues2025event}. However, these studies do not deal with constraints on actuators.

It is known that due to design or physical limitations, actuators can present restrictions in their operating range, usually modeled in the form of saturation~\cite{peixoto2022static,RefArtPedro}. If the presence of saturation is not properly addressed, the performance of the closed-loop system may deteriorate, or even lead to instability~\cite{RefArtPedro}. Therefore, the actuators' constraints need to be taken into account in the analysis and design of the control system. Motivated by this issue, some studies addressed the presence of actuator saturation in extremum seeking control systems~\cite{SaturatedES, ESinputConst}. However, the articles dealing with extremum seeking control with saturating actuators have two main limitations: $(i)$ they do not deal with the multivariable case, and $(ii)$ they do not propose control design conditions that ensure the stability of the closed-loop system even in the presence of saturation.

This article addresses the multivariable extremum seeking control problem with saturating actuators. First, for the case where the Hessian matrix is unknown, a stability analysis condition is proposed for the average system, and the convergence of the trajectories to the optimal point is ensured by invoking the averaging theorem~\cite{plotnikov1979averaging}. 
For the case where the Hessian matrix is polytopic uncertain, a control design condition is established to obtain the control gain such that the average system is exponentially stable.
Interestingly, the design methodology provided here allows the possibility of non-diagonal control gains, offering greater design flexibility than the diagonal gains usually assumed \textit{a priori} in the extremum-seeking literature.

\textbf{Notation.}
$\mathbb{R}^n$ denotes the Euclidean space $n$-dimensional and $\mathbb{R}^{m\times n}$ the set of real matrices $m\times n$. 
$X > 0~(X < 0)$ denotes that $X$ is a symmetric positive definite matrix (negative). 
For a matrix $X$, $X_{(\ell)}$ denotes its $\ell$-th row.


\section{Problem Formulation}
\label{sec:preliminaries}

\begin{figure}[!ht]
\begin{center}
\includegraphics[width=8.4cm]{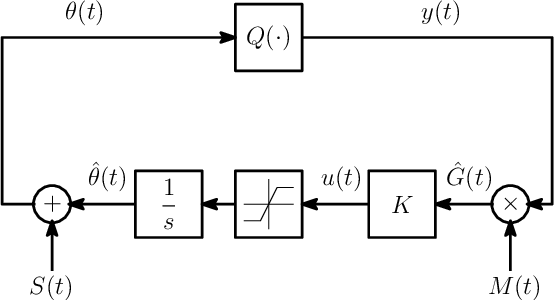}
\caption{Extremum seeking with saturated control system.} 
\label{fig:diag_ES_sat}
\end{center}
\end{figure}

Consider the extremum seeking control system with saturation based on the gradient algorithm shown in Figure~\ref{fig:diag_ES_sat}.

In this system, the multivariable quadratic static map is given by
\begin{equation} \label{eq:map_static}
       y(t) = Q(\theta(t)) =  Q^* + \frac{1}{2}(\theta(t) - \theta^*)^\top H(\theta(t) - \theta^*),
\end{equation}
where $Q^* \in \mathbb{R}^n$ is the unknown optimal point of the map, 
$\theta^* \in \mathbb{R}^n$ is the unknown optimizer of the map, 
$\theta \in \mathbb{R}^n$ is the input vector, 
$H \in \mathbb{R}^{n \times n}$ is the unknown Hessian Matrix and 
$y \in \mathbb{R}^n$ is the map output.
Even though the Hessian matrix is unknown, it can be assumed that it is a positive definite matrix when the minimum point is desired and a negative definite matrix when the maximum point is desired.

In this scheme, the signal $\theta(t)$ applied to the multivariable static map is
\begin{equation} \label{eq:input_map}
    \theta(t) = \hat{\theta}(t) + S(t)
\end{equation}
where $\hat{\theta} \in \mathbb{R}^n$ is the estimated value of $\theta^*$.
Moreover, the $\hat{\theta}$-dynamics is described as follows:
\begin{align}
    \dot{\hat{\theta}}(t) = \mathrm{sat}(u(t)) = \mathrm{sat}(K\hat{G}(t)),
\end{align}
where $K$ is the control gain, $\mathrm{sat}(\cdot)$ is the saturating function defined in the element-wise sense, as follows:
\begin{align}
    \mathrm{sat}(u) = \begin{bmatrix}
        \mathrm{sat}(u_1) \\ \vdots \\ \mathrm{sat}(u_n)
    \end{bmatrix} = 
    \begin{bmatrix}
        \mathrm{sign}(u_1) \min(|u_1|,\overline{u}_1) \\ \vdots \\ \mathrm{sign}(u_n) \min(|u_n|,\overline{u}_n)
    \end{bmatrix},
\end{align}
where $\overline{u}_\ell > 0$ is the limit of $\ell$-th control input signal, and
$\hat{G}(t)$ is the gradient estimate given by
\begin{equation}
    \hat{G}(t) = M(t)y(t).
\end{equation}
The perturbation signals of the extremum seeking are defined as follows \cite{LivroTiagoES}:
\begin{equation} \label{eq:vetor_S}
    S(t) = \begin{bmatrix}
       a_1\sin{(\omega_1 t)} & \cdots & a_n\sin{(\omega_n t)} \\
    \end{bmatrix}^\top
\end{equation}
\begin{equation} \label{eq:vetor_M}
    M(t) = \begin{bmatrix}
       \tfrac{2}{a_1}\sin{(\omega_1 t)} & \cdots & \tfrac{2}{a_n}\sin{(\omega_n t)} \\
    \end{bmatrix}^\top
\end{equation}
where $a_i$, $i = 1,\ldots,n$, are the non-zero amplitudes, and the frequencies of the disturbance signals are selected such that 
\begin{align} \label{eq:relacao_sinais_omega_pert}
    \omega_i = \omega_i' \omega, \quad i = 1,\ldots,n,
\end{align}
e $\omega_i' \notin \{\omega_j', \frac{1}{2}(\omega_j'+\omega_k'), \omega_k'\pm \omega_l' \}$,
for all $i,j,k = 1,\ldots,n$.

Defining the estimation error
\begin{equation} \label{eq:error_map}
    \tilde{\theta}(t) = \hat{\theta}(t) - \theta^*,
\end{equation}
we can compute its dynamics as
\begin{align} \label{eq:theta_til_ponto_&_saturacao}
    \dot{\tilde{\theta}}(t) = \dot{\hat{\theta}}(t) = \mathrm{sat}(u(t)) = \mathrm{sat}(K\hat{G}(t)).
\end{align}
Thus, it is clear that the estimation error $\tilde{\theta}(t)$ converges towards zero if the gradient estimate $\hat{G}(t)$ converges to zero.
\textcolor{black}{In contrast to previous work, this article considers the presence of saturation in the signal $u(t)$, which can be viewed as an ``input signal'' to the extremum-seeking control system~\cite{rodrigues2025event}. 
Note also that dealing with the saturation before the integration can be
a more challenging problem than considering the saturation directly in the input of the static map, 
as done in~\cite{ESinputConst}. The reason is that the saturation before the integration might lead to the
well-known windup phenomenon, which can deteriorate the closed-loop performance.}

As $y(t)$ is the quadratic map output, the gradient estimate can be written as:
\begin{equation}
    \hat{G}(t) = M(t)\left(Q^* + \frac{1}{2}(\tilde{\theta}(t) + S(t))^\top H(\tilde{\theta}(t) + S(t))\right),
\end{equation}
or
\begin{multline} \label{eq:grad_simp}
     \hat{G}(t) = M(t)Q^* + \frac{1}{2}M(t)\tilde{\theta}^\top(t)H\tilde{\theta}(t)\\
     + M(t)S^\top(t)H\tilde{\theta}(t) + \frac{1}{2}M(t)S^\top(t)HS(t).
\end{multline}
By defining the matrix
\begin{equation} \label{eq:delta}
    \Omega(t) = M(t)S^\top(t)H,
\end{equation}
the multiplication in \eqref{eq:delta} results in a matrix of the following form:
\begin{equation} \label{eq:delta_final}
    \Omega(t) = H + \Delta(t)H,
\end{equation}
where
    $\Delta_{ii} = 1-\cos(2\omega_it)$,
    $\Delta_{ij} = \frac{a_{j}}{a_{i}}\cos(\omega_i-\omega_j) - \frac{a_{j}}{a_{i}}\cos(\omega_i+\omega_j)$.

Due that, \eqref{eq:grad_simp} can be expressed as
\begin{multline} \label{eq:grad_simp_2}
     \hat{G}(t) = M(t)Q^* + \frac{1}{2}M(t)\tilde{\theta}^\top H\tilde{\theta}(t) \\
     +  \Omega(t)\tilde{\theta}(t) + \frac{1}{2}\Omega(t)S(t).
\end{multline}
Then, the dynamics of \eqref{eq:grad_simp_2} can be rewritten as:
\begin{equation} \label{eq:grad_derivado}
     \dot{\hat{G}}(t) = H\mathrm{sat}(u(t)) + \Delta(t)H\mathrm{sat}(u(t)) + \varsigma(t)
\end{equation}
where
\begin{multline} \label{eq:resto_grad_derivado}
        \varsigma(t) = \dot{M}(t)Q^* + \dot{\Delta}(t)H\tilde{\theta}(t) + \frac{1}{2}H\dot{S}(t)+ \frac{1}{2} \dot{\Delta}(t)HS(t)  \\ + \frac{1}{2}\Delta(t)H\dot{S}(t) + [\frac{1}{2}M(t)\tilde{\theta}^\top H\tilde{\theta}(t)].
\end{multline} 

\subsection{Defining a new time-scale}

For the analysis of the closed-loop system stability, a change in the time scale is performed. From the relation of the disturbance signal frequencies \eqref{eq:relacao_sinais_omega_pert}, it is ensured that the frequency ratio must be rational. Thus, there exists a period $T$ that
\begin{align}
    T = 2\pi \times \mathrm{LCM}\left\lbrace\frac{1}{\omega_i}\right\rbrace, \quad  i  = 1, 2, \ldots,n,
\end{align}
where $\mathrm{LCM}$ denotes the least common multiple. The change of time scale of the system in \eqref{eq:grad_derivado} consists of a transformation $\tau = \omega t$, where
\begin{align}
    \omega := \frac{2\pi}{T}.
\end{align}

So, the system \eqref{eq:grad_derivado} can be rewritten as
\begin{align} \label{eq:esc_temp_mudanca_Grad}
    \frac{d\hat{G}\left(\tau\right)}{d\tau} = \frac{1}{\omega}\mathcal{F}\left(\tau, \hat{G}, \tilde{\theta}, \frac{1}{\omega}\right) 
\end{align}
where 
\begin{align}\label{eq:esc_temp_mudanca_Grad_function}
  \mathcal{F}\left(\tau,  \hat{G}, \tilde{\theta}, \frac{1}{\omega}\right)   = H\mathrm{sat}(u(\tau)) + \Delta(\tau)H\mathrm{sat}(u(\tau)) + \varsigma(\tau)
\end{align}

\subsection{Average System}

Calculating the average version of \eqref{eq:esc_temp_mudanca_Grad}, we have 
\begin{align} \label{eq:Sistema_Medio_Derivado}
    \frac{d\hat{G}_\mathrm{av}(\tau)}{d\tau} = \frac{1}{\omega}\mathcal{F}_\mathrm{av}(\hat{G}_\mathrm{av}) 
\end{align}
where
\begin{align} \label{eq:Media_calculada_Sistema}
    \mathcal{F}_\mathrm{av}(\hat{G}_\mathrm{av}) = \frac{1}{T} \int_0^T   \mathcal{F}_\mathrm{av}(\delta, \hat{G}_\mathrm{av}, 0) d\delta.
\end{align}
For each term, the average is computed below:
\begin{align} 
    S_\mathrm{av}(\tau) &= \frac{1}{T} \int_0^T   S(\delta) d\delta = 0, \dot{S}_\mathrm{av}(\tau) = \frac{1}{T} \int_0^T   \dot{S}(\delta) d\delta = 0, \label{eq:Media_Sinal_S} \\
   M_\mathrm{av}(\tau) &= \frac{1}{T} \int_0^T   M(\delta) d\delta = 0, \dot{M}_\mathrm{av}(\tau) = \frac{1}{T} \int_0^T   \dot{M}(\delta) d\delta = 0, \label{eq:Media_Sinal_M} \\
    \Delta_\mathrm{av}(\tau) &= \frac{1}{T} \int_0^T   \Delta(\delta) d\delta = 0, \dot{\Delta}_\mathrm{av}(\tau) = \frac{1}{T} \int_0^T   \dot{\Delta}(\delta) d\delta = 0. \label{eq:Media_Sinal_Delta}
\end{align}
As a result, one can obtain
\begin{align*}
    \Omega_{\mathrm{av}}(\tau) = \frac{1}{T} \int_0^T \Omega(\delta) d\delta = H,  
    \dot{\Omega}_{\mathrm{av}}(\tau) = \frac{1}{T} \int_0^T \dot{\Omega}(\delta) d\delta = 0.
\end{align*}
Then, the average system is finally given by
\begin{equation} \label{eq:grad_medio}
     \dot{\hat{G}}_{\mathrm{av}}(\tau) = \frac{1}{\omega} H\mathrm{sat}(u_{\mathrm{av}}(\tau)) = \frac{1}{\omega}H\mathrm{sat}(K\hat{G}_{\mathrm{av}}(\tau)).
\end{equation}

Consider the dead zone nonlinearity of the control input signal~\cite{LivroSaturacao}:
\begin{equation} \label{eq:def_sat}
    \psi(u) = u - \mathrm{sat}(u).
\end{equation}
Based on this, the average closed-loop system obtained from~\eqref{eq:grad_medio} and \eqref{eq:def_sat} can be written as
\begin{align}\label{eq:average_system1}
    \dot{\hat{G}}_{\mathrm{av}}(\tau) = \frac{1}{\omega}H K \hat{G}_{\mathrm{av}}(\tau) - \frac{1}{\omega}H \psi(u_{\mathrm{av}}(\tau)),
\end{align}
where $u_{\mathrm{av}} = K \hat{G}_{\mathrm{av}}$.

In general, solutions available in the literature are developed for the stability analysis of extremum seeking control systems, assuming the knowledge of the sign of the Hessian matrix $H$. Based on this, a diagonal structure with the opposite sign is assigned to the gain matrix $K$. Although this approach requires little knowledge of the Hessian matrix $H$, it becomes difficult to design the gain matrix using constructive design techniques via LMIs. 

For this purpose, it is assumed that the Hessian matrix $H$ is unknown, but takes values within a polytopic set according to the following parameterization:
\begin{align}\label{eq:polytpoic_H}
    H = H(\alpha) = \sum_{i=1}^N \alpha_i H_i,
\end{align}
where the vector of uncertain parameters $\alpha = (\alpha_1,\ldots,\alpha_N)$ belongs to the unitary simplex
\begin{align}
    \Lambda = \left\lbrace \alpha \in \mathbb{R}^N : \sum_{i = 1}^N \alpha_i = 1, \; \alpha_i \geq 0, i = 1,\ldots,N \right\rbrace
\end{align}
and $H_i \in \mathbb{R}^{n \times n}$, $i = 1,\ldots, N$ are the polytope vertices, that are known matrices. Thus, it is possible to obtain the following uncertain polytopic description for the closed-loop system:
\begin{align} \label{eq:dinamica_Gav_politopo}
    \dot{\hat{G}}_{\mathrm{av}}(\tau) = \frac{1}{\omega}H(\alpha) K \hat{G}_{\mathrm{av}}(\tau) - \frac{1}{\omega}H(\alpha) \psi(u_{\mathrm{av}}(\tau)).
\end{align}

The problem addressed in this paper is to design a control robust gain $K \in \mathbb{R}^{n \times n}$ such that the closed-loop average system is exponentially stable. Then, by invoking the averaging theorem, we prove the stability of the non-averaged closed-loop system~\eqref{eq:grad_derivado}.


\section{Main Results}
\label{sec:Main_Results}

This section presents a stability condition for the average system based on a sector condition for the dead-zone function. After that, a convergence condition of the original system is established using the averaging theorem of~\cite{plotnikov1979averaging}. Finally, the control design condition is formulated as an optimization problem based on LMIs.

\subsection{Stability Analysis of the Average System}

The following lemma establishes a sector condition for the dead-zone nonlinearity.
\begin{lema}\label{lem:saturacao}
    Consider a diagonal positive definite matrix $U \in \mathbb{R}^{n \times n}$ and a matrix $L \in \mathbb{R}^{m \times n}$. If $\hat{G}_{\mathrm{av}}$ belongs to the set
    \begin{align}\label{eq:conjunto_G}
        \mathcal{G} = \left\lbrace \hat{G}_{\mathrm{av}} \in \mathbb{R}^n : |(K-L)_{(\ell)} \hat{G}_{\mathrm{av}}| \leq \overline{u}_{\ell}, \, \ell = 1,\ldots,m \right\rbrace,
    \end{align}
    then
    \begin{align}\label{eq:setor_psi}
        \psi^\top(u_{\mathrm{av}}) U \left(\psi(u_{\mathrm{av}}) - L \hat{G}_{\mathrm{av}}\right) \leq 0.
    \end{align}
\end{lema}
\begin{proof}
    The demonstration follows similar steps to~\cite[Lemma~1]{RefArtPedro}.
\end{proof}

\begin{teorema}\label{thm:1}
    Consider the average system in~\eqref{eq:average_system1} of the extremum-seeking control system subject to saturation. 
    Given a positive scalar $\eta > 0$, a control gain $K \in \mathbb{R}^{n \times n}$ and matrices $P = P^\top > 0 \in \mathbb{R}^{n \times n}$,
    $Q = Q^\top > 0 \in \mathbb{R}^{n \times n}$, there are matrices
    $L \in \mathbb{R}^{n \times n}$, and a diagonal matrix $U > 0 \in \mathbb{R}^{n \times n}$ such that
    \begin{align}\label{eq:ineq-stability}
        \begin{bmatrix}
            K^\top H P + P H K + 2 \eta P & L^\top U - P H \\
            U L - H P & - 2 U
        \end{bmatrix} < 0,
    \end{align}
    then the origin of the average system is exponentially stable.
\end{teorema}
\begin{proof}
Assume that the condition~\eqref{eq:ineq-stability} is satisfied.
Multiplying by $[\hat{G}_{\mathrm{av}}^\top \; \psi^\top(u_\mathrm{av})]$ on the left and by its transpose on the right results in
\begin{align}
     &\hat{G}_{\mathrm{av}}^\top \left(PHK + K^\top H P\right)\hat{G}_{\mathrm{av}} - 2 \hat{G}_{\mathrm{av}}^\top P H \psi(u_\mathrm{av}) \nonumber \\
     &- 2 \psi^\top(u_{\mathrm{av}}) U \left(\psi(u_{\mathrm{av}}) - L \hat{G}_{\mathrm{av}}\right) + 2 \eta \hat{G}_{\mathrm{av}}^\top P \hat{G}_{\mathrm{av}} < 0.
\end{align}
Thus, it follows from Lemma~\ref{lem:saturacao} that if $\hat{G}_{\mathrm{av}} \in \mathcal{G}$, with $\mathcal{G}$ given by~\eqref{eq:conjunto_G}, then 
\begin{align}\label{eq:exp-ineq}
    \dot{V}(\hat{G}_{\mathrm{av}}) \leq - 2 \eta V(\hat{G}_{\mathrm{av}}) < 0,
\end{align}
where
\begin{align} \label{eq:lyapunov_Gav}
    V(\hat{G}_{\mathrm{av}}) = \hat{G}_{\mathrm{av}}^\top P \hat{G}_{\mathrm{av}},
\end{align}
is a Lyapunov function that ensures the exponential stability of the origin of the average system.
From the Comparison Lemma, it follows from~\eqref{eq:exp-ineq} that
\begin{align}
    {V}(\hat{G}_{\mathrm{av}}(\tau)) \leq e^{-2 \eta t} {V}(\hat{G}_{\mathrm{av}}(0))
\end{align}

Furthermore, as 
\begin{align}
    \lambda_{\min}(P) \|\hat{G}_{\mathrm{av}}\|^2 \leq 
    V(\hat{G}_{\mathrm{av}}) \leq \lambda_{\max}(P) \|\hat{G}_{\mathrm{av}}\|^2,
\end{align}
it can be obtained
\begin{align}
    \|\hat{G}_{\mathrm{av}}\| \leq \kappa e^{-\eta \tau}\|\hat{G}_{\mathrm{av}}(0)\|
\end{align}
where $\kappa = \sqrt{\lambda_{\max}(P) / \lambda_{\min}(P)}$.
Then, the origin of the system is exponentially stable.
\end{proof}
\begin{lema}
 Let
 \begin{equation}
   V(\hat{G}_{\mathrm{av}}) = \hat{G}_{\mathrm{av}}^{\top} P \hat{G}_{\mathrm{av}},
 \end{equation}
and consider the following set
\begin{align}\label{eq:elipse_V}
   \mathcal{E} = \{\hat{G}_{\mathrm{av}} \in \mathbb{R}^n : \hat{G}_{\mathrm{av}}^{\top} P \hat{G}_{\mathrm{av}} \leq 1 \}.
\end{align}
If the following conditions are satisfied
 \begin{equation}\label{eq:sat}
 \begin{bmatrix}
   P     &  K_{(\ell)}^\top - L_{(\ell)}^\top \\
   \star & \bar u_{(\ell)}^{2}                                 
 \end{bmatrix} >  0, \quad \ell  = 1,2,\ldots,n,
\end{equation} 
and $\hat{G}_{\mathrm{av}}(0)$ is taken in the region
$\mathcal{E}$ in~\eqref{eq:elipse_V}, then $\mathcal{E} \subset \mathcal{G}$, with $\mathcal{G}$ in~\eqref{eq:conjunto_G}
\end{lema}
\begin{proof} Consider \( V(\hat{G}_{\mathrm{av}}) = \hat{G}_{\mathrm{av}}^{\top} P \hat{G}_{\mathrm{av}} \) and the following relation:
\begin{equation}\label{eq:sat2_condition}
 V(\hat{G}_{\mathrm{av}}) \geq \frac{|(K-L)_{(\ell)} \hat{G}_{\mathrm{av}}|^2}{\overline{u}_{(\ell)}^2}.
\end{equation}
Thus, the relation in \eqref{eq:sat2_condition} can be written as:
\begin{equation}
 \hat{G}_{\mathrm{av}}^{\top}(\tau) P \hat{G}_{\mathrm{av}} \geq \overline{u}_{(\ell)}^{-2} ((K-L)_{(\ell)} \hat{G}_{\mathrm{av}})^\top (K-L)_{(\ell)} \hat{G}_{\mathrm{av}},
\end{equation}
or still
\begin{equation}
 \hat{G}_{\mathrm{av}}^{\top}(\tau) \left( P  - \overline{u}_{(\ell)}^{-2} (K_{(\ell)}-L_{(\ell)})^\top (K_{(\ell)}-L_{(\ell)}) \right)\hat{G}_{\mathrm{av}} \geq 0.
\end{equation}
Applying the Schur complement results in~\eqref{eq:sat}.
This means that if $\hat{G}_{\mathrm{av}}$ is taken in the region
$\mathcal{E}$ em~\eqref{eq:elipse_V}, so $\mathcal{E} \subset \mathcal{G}$. This concludes the proof.
\end{proof}

\subsection{Stability Analysis Using Averaging Theorem}

\begin{teorema}
    Consider the average closed-loop dynamic of the gradient estimate subject to saturation~\eqref{eq:average_system1}.
    If the theorem conditions of Theorem~\ref{thm:1} are satisfied, 
    then, for $\omega > 0$ sufficiently large, the equilibrium $\hat{G}_{\mathrm{av}} = 0$ is exponentially stable and $\tilde{\theta}_{\mathrm{av}}(t)$ converges exponentially to zero. In particular, there exist constants $\overline{\kappa}, \kappa_y > 0$ such as
    \begin{align} \label{eq:theta_desigualdades}
        \|\theta(t) - \theta^\ast\| \leq \overline{\kappa} e^{-\eta t} + \mathcal{O}\left(a + \frac{1}{\omega}\right) \\
        \label{eq:y_desigualdades}
        |y(t) - Q^\ast| \leq  \kappa_y e^{-\eta t} + \mathcal{O}\left(a^2 + \frac{1}{\omega^2}\right),
    \end{align}
    where $a = \sqrt{\sum_{i=1}^n a_i^2}$, taking $a_i$ the defined constants in~\eqref{eq:vetor_S} and $\overline{\kappa}$ and $\kappa_y$ constants which depends on the initial condition $\theta(0)$. 
\end{teorema}
\begin{proof}
    From the equation~\eqref{eq:grad_simp_2}, it can be obtained that
    \begin{align}
        \hat{G}_{\mathrm{av}}(\tau) = \frac{1}{\omega}H \tilde{\theta}(\tau),
    \end{align}
    since the quadratic term $\frac{1}{2}M(t)\tilde{\theta}^\top H\tilde{\theta}(t)$ can be neglected in a local analysis, and the other terms have zero average.
    
    Rewriting the Lyapunov function in~\eqref{eq:lyapunov_Gav} as
    \begin{align}
        V(\tilde{\theta}_{\mathrm{av}}) = \tilde{\theta}_{\mathrm{av}}^\top \overline{P} \tilde{\theta}_{\mathrm{av}},
    \end{align}
    where $\overline{P} = H^\top P H$ is an symmetric positive definite matrix, provided that $H$ and $P$ are symmetric and positive definite. Thus, it is possible to find
    \begin{align}
        \|\tilde{\theta}_{\mathrm{av}}(\tau)\| \leq \overline{\kappa} e^{-\eta \tau / \omega} \|\tilde{\theta}_{\mathrm{av}}(0)\|,
    \end{align}
    where $\overline{\kappa} = \sqrt{\lambda_{\max}(\overline{P}) / \lambda_{\min}(\overline{P})}$. As the differential equation in~\eqref{eq:esc_temp_mudanca_Grad} has discontinuity on the right side, due to the presence of the saturating function, \eqref{eq:esc_temp_mudanca_Grad_function} is $T$-periodic and Lipschitz continuous, it follows from~\cite{plotnikov1979averaging}
    that 
    \begin{align}
        \|\tilde{\theta}(t) - \tilde{\theta}_{\mathrm{av}}(t)\| \leq \mathcal{O}\left(\frac{1}{\omega}\right).
    \end{align}

    Applying the triangular inequality, it can be guaranteed that
    \begin{align} 
        \| \tilde{\theta}(t) \| 
        &\leq \overline{\kappa} e^{-\eta t} \| \tilde{\theta}_{\text{av}}(0) \| + \mathcal{O}\left(\frac{1}{\omega}\right).
    \end{align}

    Applying the averaging theorem~\cite{plotnikov1979averaging}, it can be concluded that
    \begin{align}
        \| \hat{G}(t) - \hat{G}_{\text{av}}(t) \| \leq \mathcal{O}\left(\frac{1}{\omega}\right).
    \end{align}
    Similarly, we can apply the triangular inequality to obtain
    \begin{align}
        \| \hat{G}(t) \| 
        &\leq \overline{\kappa} e^{-\eta t} \| \hat{G}_{\text{av}}(0) \| + \mathcal{O}\left(\frac{1}{\omega}\right).
    \end{align}
    From~\eqref{eq:input_map} and the definition of $\tilde{\theta}(t)$, we have
    \begin{align}
        \theta(t) - \theta^* = \tilde{\theta}(t) + S(t).
    \end{align}

Thus, the following relation can be obtained:
    \begin{align} \label{eq:demonstracao_estabilidade_theta}
        \| \theta(t) - \theta^* \| \leq (\overline{\kappa}) e^{-\eta t} \|\theta(0)- \theta^*\| + \mathcal{O}\left(a + \frac{1}{\omega}\right)
    \end{align}

Let the error variable
\begin{align}
    \tilde{y}(t) := y(t) - Q^*,~~~~y(t) = Q(\theta(t)).
\end{align}
By computing its norm, and using the Cauchy–Schwarz inequality, one gets
\begin{align}
    |\tilde{y}(t)| & = |y(t) - Q^*| = |(\theta(t) - \theta^*)^\top H(\theta(t) - \theta^*)| \\ &\leq \|H\| \|((\theta(t) - \theta^*))\|^2.
\end{align}
From~\eqref{eq:demonstracao_estabilidade_theta}, it is still possible to obtain
\begin{align}
    |\tilde{y}(t)| \leq \|H\| ((\overline{\kappa})^2 e^{-2\eta t} \|\theta(0)- \theta^*\|^2 + \mathcal{O}\left(a^2 + \frac{2a}{\omega} + \frac{1}{\omega^2}\right)
\end{align}

As $e^{-\eta t} \geq e^{-2 \eta t}$ for $\omega > 0$, and $a^2+\frac{1}{\omega^2} \geq \frac{2a}{\omega}$, by the Young's inequality, one obtains
\begin{align}
    |y(t) - Q^*| \leq \kappa_y e^{-\eta t} + \mathcal{O}\left(a^2 + \frac{1}{\omega^2}\right),
\end{align}
where 
\begin{align*}
    \kappa_y = \|H\| (\overline{\kappa})^2  \|\theta(0)- \theta^*\|^2
\end{align*}

As a result, the inequalities \eqref{eq:theta_desigualdades} and \eqref{eq:y_desigualdades} are guaranteed. This concludes the proof.
\end{proof}

\subsection{Control Design Condition}

The theorem below provides a constructive LMI-based condition for designing the gain of the extremum-seeking control system.
\begin{teorema}\label{thm:2}
Let $\eta > 0$ be a given scalar.
If there exist a symmetric positive definite matrix $W \in \mathbb{R}^{n \times n}$, diagonal positive definite matrix $V \in \mathbb{R}^{n \times n}$, and matrices $Z, Y, T \in \mathbb{R}^{n \times n}$, such that the inequalities below are satisfied for all $1,2,\ldots,N$:
\begin{equation}\label{eq:LMI_thm}
    \begin{bmatrix} 
        H_i Z+Z^\top H_i ^\top + 2\eta W & \star & \star\\
         W-T^{\top}+\epsilon H_i Z & -\epsilon (T^{\top}+T) & \star \\
         Z + Y - VH_i ^\top & -\epsilon V H_i ^\top & -2V
    \end{bmatrix}<0,
\end{equation}
\begin{align}
    \begin{bmatrix}
            W &  Z_{(\ell)}^\top - Y_{(\ell)}^\top \\
            \star & u_\ell^2
    \end{bmatrix} \geq 0, \quad \ell = 1,\ldots,n,  
    \label{eq:saturation_thm} 
\end{align}
Therefore, if the conditions in~\eqref{eq:LMI_thm} are satisfied, so the origin of the average uncertain closed-loop system~\eqref{eq:dinamica_Gav_politopo} with $K = Z T^{-1}$ is exponentially stable.
\end{teorema}

\begin{proof}
Take the following candidate Lyapunov function:
\begin{align} \label{eq:lyapunov_Gav_2}
    V(\hat{G}_{\mathrm{av}}) = \hat{G}_{\mathrm{av}}^\top P \hat{G}_{\mathrm{av}},
\end{align}
which is positive definite for all $\hat{G}_{\mathrm{av}} \neq 0 \in \mathbb{R}^n$, for an matrix $P = P^\top > 0$.

Calculating the temporal derivative of \eqref{eq:lyapunov_Gav} and applying the S-procedure with~\eqref{eq:setor_psi}, there is

\begin{align} \label{eq:dev_lyapunov_Gav}
    \dot{V}(\hat{G}_{\mathrm{av}}) - 
    2 \psi^\top(u_{\mathrm{av}}) U
        \left(\psi(u_{\mathrm{av}}) -L\hat{G}_{\mathrm{av}}\right) + 2 \eta V(\hat{G}_{\mathrm{av}}) < 0,
\end{align}
that is,
\begin{multline}\label{eq:V_diff}   \dot{\hat{G}}_{\mathrm{av}}^\top P \hat{G}_{\mathrm{av}} + \hat{G}_{\mathrm{av}}^\top P \dot{\hat{G}}_{\mathrm{av}} - \\
       2 \psi^\top(u_{\mathrm{av}}) U
        \left(\psi(u_{\mathrm{av}}) -L\hat{G}_{\mathrm{av}}\right) +  2\eta \hat{G}_{\mathrm{av}}^\top P \hat{G}_{\mathrm{av}}< 0.
\end{multline}

Consider the augmented vector shown below:
\begin{align}
    \xi = 
    \begin{bmatrix}
        \hat{G}_{\mathrm{av}} \\ \dot{\hat{G}}_{\mathrm{av}} \\
        \psi(u_{\mathrm{av}})
    \end{bmatrix}.
\end{align}

It follows from~\eqref{eq:V_diff} that
\begin{align}
    \xi^\top 
    \begin{bmatrix}
        2\eta P & P & L^\top U\\
        P & 0 & 0\\
        UL & 0 & -2U
    \end{bmatrix}\xi < 0
\end{align}
for all $\xi \neq 0$ and $\mathcal{B} \xi = 0$, where
\begin{align}
    \mathcal{B} = 
    \begin{bmatrix}
        H(\alpha)K & -I & -H(\alpha)
    \end{bmatrix}.
\end{align}
Using the Finsler Lemma, we have that
\begin{align*}
    \begin{bmatrix}
        2\eta P & P & L^\top U\\
        P & 0 & 0\\
        UL & 0 & -2U
    \end{bmatrix} + 
    \begin{bmatrix}
        X^\top \\ \epsilon X^\top \\ 0
    \end{bmatrix}
     \begin{bmatrix}
        H(\alpha)K & -I & -H(\alpha)
    \end{bmatrix}
     \\
    + \begin{bmatrix}
        K^\top H^\top(\alpha) \\ -I \\ -H^\top(\alpha)
    \end{bmatrix}
    \begin{bmatrix}
        X & \epsilon X & 0
    \end{bmatrix} < 0
\end{align*}
 By applying a congruence transformation with
 \begin{align}
     \begin{bmatrix}
         X^{-\top} & 0 & 0 \\
         0 & X^{-\top} & 0 \\
         0 & 0 & U^{-1}
     \end{bmatrix}
 \end{align}
and using the transformations of variables $Z = KX^{-1}$, $W=X^{-\top}PX^{-1}$, $T = X^{-1}$, $V = U^{-1}$ and $Y = LX^{-1}$, the resulting inequality is
\begin{align*}
    \begin{bmatrix}
        HZ+Z^\top H^\top + 2\eta W &\star &  \star \\
         W-T^{\top}+\epsilon H Z & -\epsilon ( T^{\top}+ T) & \star \\
         Y - VH^\top & -\epsilon V H^\top & -2V
    \end{bmatrix}<0.
\end{align*}

Multiplying \eqref{eq:saturation_thm} by~$\mathrm{diag}(X^{\top},1)$ on the left and by its transpose on the right, results in~\eqref{eq:sat}. This guarantees that~$\mathcal{E} \subset \mathcal{G}$. 
This concludes the proof.
\end{proof}

\section{Numerical Results}

Consider the extremum-seeking control system with non-linear map \eqref{eq:map_static} with unknown Hessian matrix taking values in the polytopic set given by the following vertices

\begin{align}
    H_1 = (1-\overline{\delta}) H_0, \quad  H_2 = (1+\overline{\delta}) H_0,
\end{align}
where $\overline{\delta} > 0$ is a parameter used to construct the vertices of the polytopic domain and $H_0$ is the Hessian matrix used in~\cite{MultVar_ESC}:
\begin{align}
    H_0 = \begin{bmatrix}
        100 & 30 \\
        30 & 20
    \end{bmatrix} > 0.
\end{align}
In addition, for the simulation, it was assumed that unknown parameters are $Q^* = 10$ and $\theta^*$ =
$\begin{bmatrix}
    2 & 4 \\
\end{bmatrix}^\top$. 
Note that the unknown parameters $Q^*$ and $\theta^*$ are not used to design the control gain, only for system simulation.

\begin{figure}[ht!]
     \centering
    \subfloat[\label{fig:sinal_controle_LMI}$\mathrm{sat}(u(t))$ -- {Theorem~\ref{thm:2}}]{
    \includegraphics[width=0.235\textwidth]{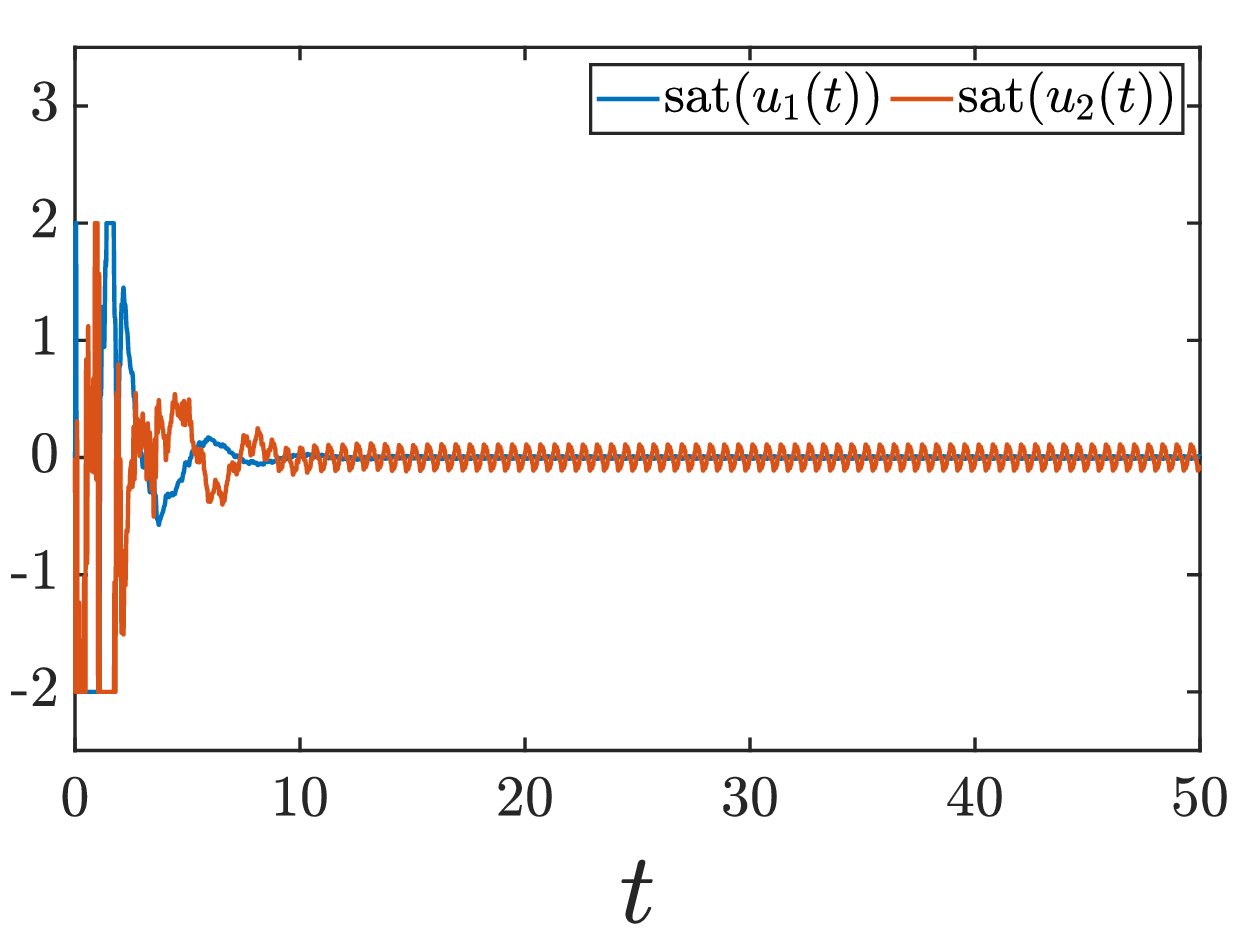}}
    \hfill
    \subfloat[\label{fig:sinal_controle_Ganho_Diagonal}$\mathrm{sat}(u(t))$ -- {Diagonal Gain}]{
    \includegraphics[width=0.235\textwidth]{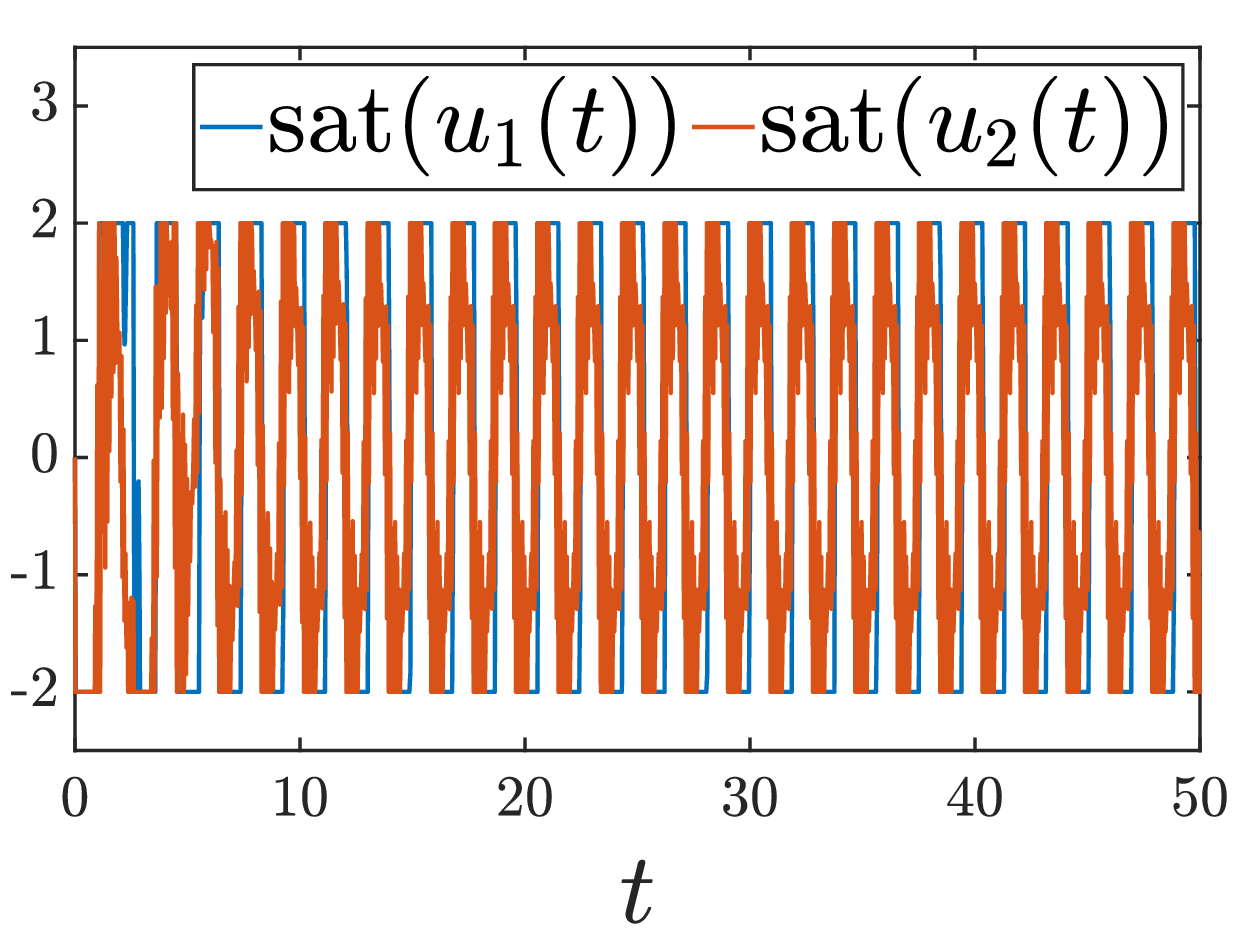}}

\vskip 3mm
    \subfloat[\label{fig:sinal_theta_LMI}$\theta(t)$ -- {Theorem~\ref{thm:2}}]{
    \includegraphics[width=0.235\textwidth]{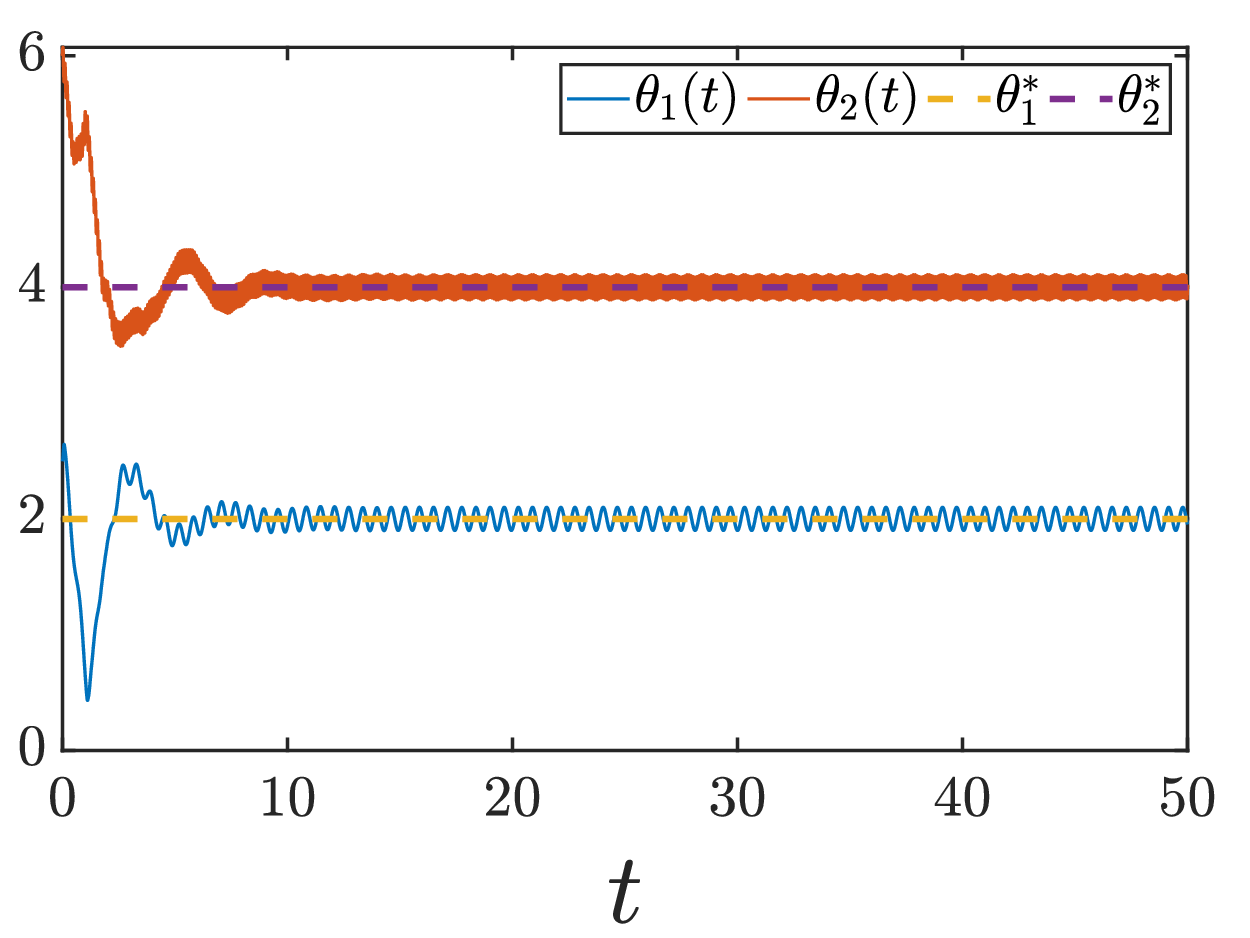}}
    \hfill
    \subfloat[\label{fig:sinal_theta_Ganho_Diagonal}$\theta(t)$ -- {Diagonal Gain}]{
    \includegraphics[width=0.235\textwidth]{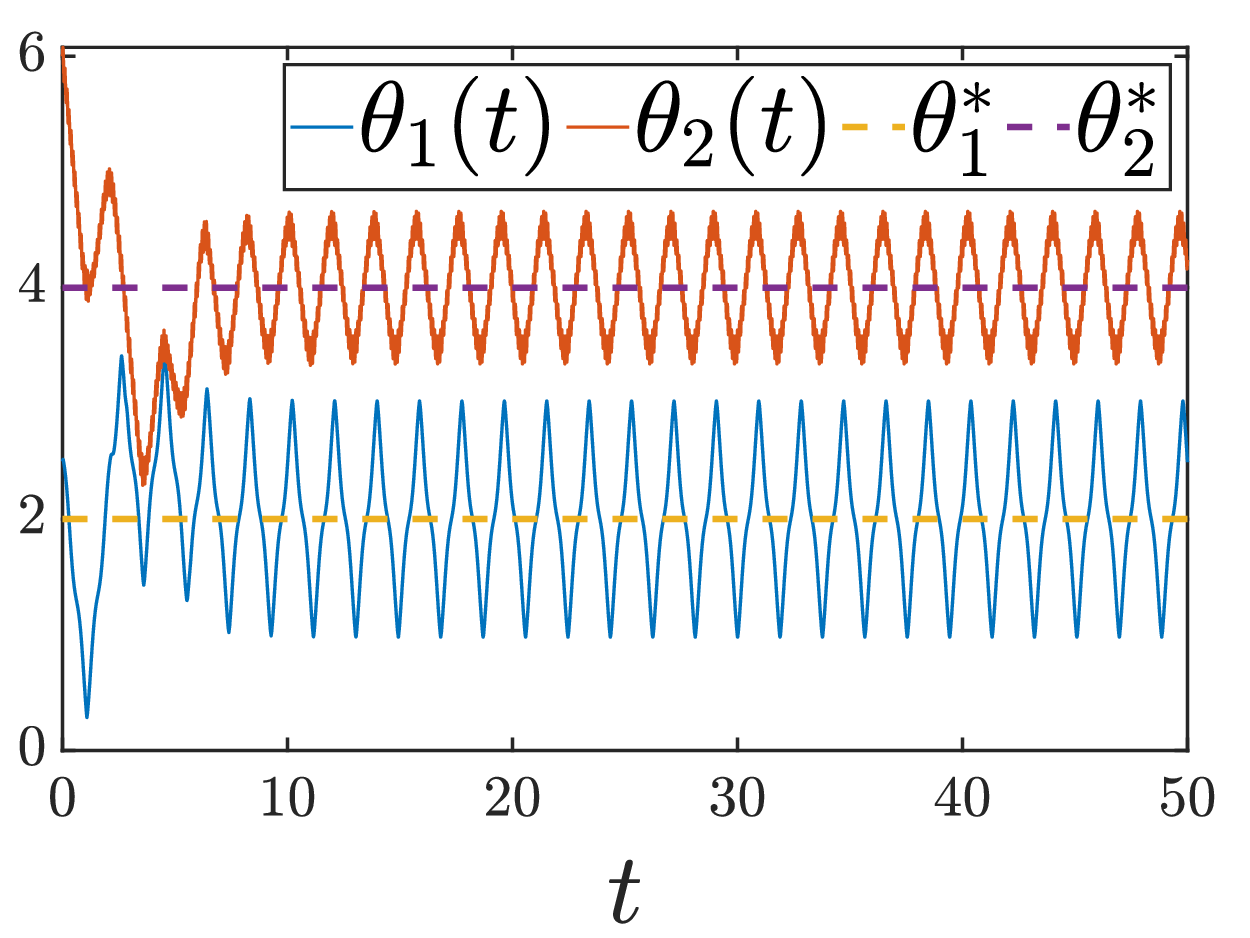}}\\

\vskip 3mm
    \subfloat[\label{fig:sinal_saida_LMI}$y(t)$ -- {Theorem~\ref{thm:2}}]{
    \includegraphics[width=0.235\textwidth]{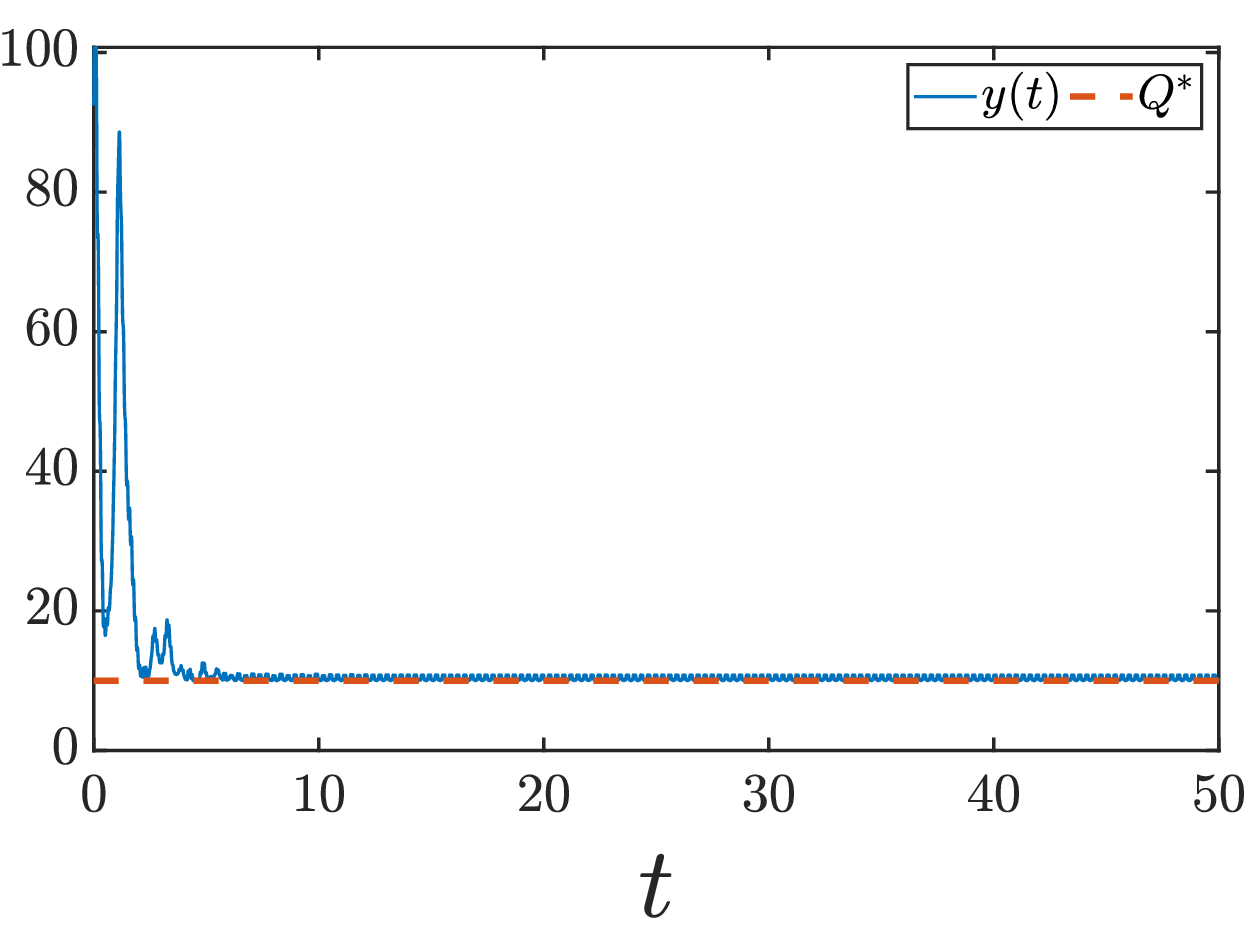}}
    \hfill
    \subfloat[\label{fig:sinal_y_Ganho_Diagonal}$y(t)$ -- {Diagonal Gain}]{
    \includegraphics[width=0.235\textwidth]{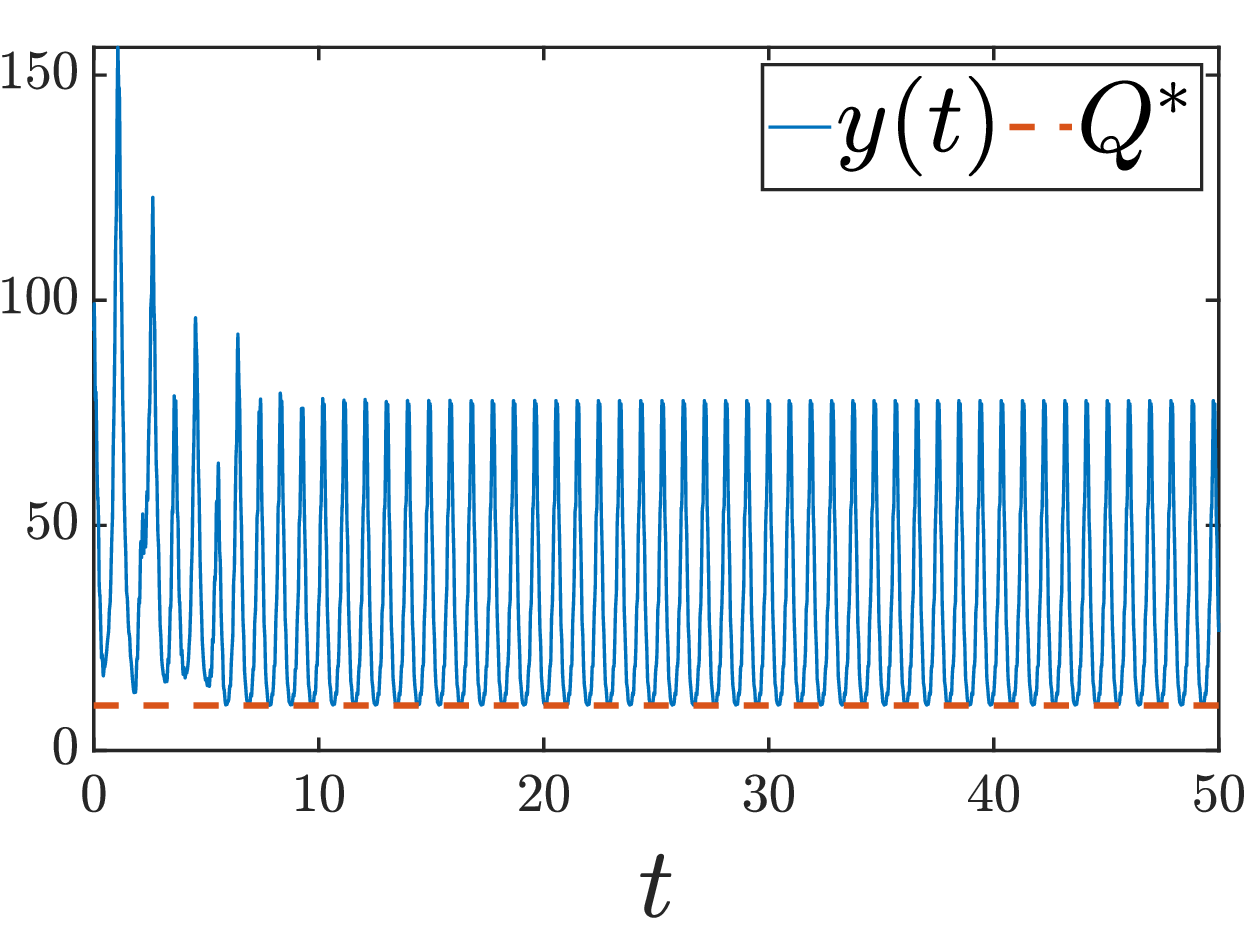}}
    
\caption{Responses of the closed-loop system with the gain designed with Theorem~\ref{thm:2} and the diagonal gain.}
\label{fig:comparacao_ganho_LMI_arbitrario_Diag}
\end{figure}

The design was performed by solving the optimization problem 
\begin{align*}
    &\max~~\mathrm{logdet}(Q_0) \\
    &\mathrm{sujeito~a}~Q_0 > 0, W \geq Q_0, \eqref{eq:LMI_thm}-\eqref{eq:saturation_thm},
\end{align*}
considering the saturation levels $\overline{u}_1 = \overline{u}_2 = 2$, $\overline{\delta} = 0.1$, $\epsilon = 0.5$ and $\eta = 1$. The controller designed was
$$K = \begin{bmatrix}
    -0.0662   & 0.0666 \\
    0.0960 &  -0.3655
\end{bmatrix}.$$
For the simulation, the frequencies of the disturbance vectors~\eqref{eq:vetor_S} and~\eqref{eq:vetor_M} are
selected as $\omega_1 = 50$ rad/s and $\omega_2 = 70~\mathrm{rad/s}$, the disturbance frequencies are $a_1=a_2=0.1$.
In the Figure~\ref{fig:comparacao_ganho_LMI_arbitrario_Diag}, the extremum-seeking control system with saturation was simulated with the designed controller using the conditions in Theorem~\ref{thm:2}, developed in this work, and with a negative diagonal gain given by $K = -0.02 I_2$. 
The choice of a diagonal gain is common in the extremum-seeking control system literature. The simulations were performed considering the initial condition $\theta(0) =  [2.5 \; 6]^\top$. As a result, it is noted that the system does not converge using the diagonal structure gain. However, with the gain designed with the proposed conditions, it was possible to ensure the convergence of the extremum-seeking control system with saturating actuators. 

\section{Conclusion}
\label{sec:conclusion}
This paper addresses the problem of multivariable extremum control subject to actuator saturation. Using a sector representation, a stability analysis condition was established for the mean system under saturation. By invoking the averaging theorem for systems with right-hand discontinuities, it was ensured that the system trajectories converge to a neighborhood of the unknown optimal point. Furthermore, assuming an uncertain polytopic representation for the Hessian matrix, a constructive condition was derived for designing the stabilizing controller. Numerical simulations illustrated the effectiveness of the proposed controller by comparing it with a diagonal gain, which is the structure commonly employed in the extremum-seeking control literature. Future research lies in the design and analysis of different control problems with saturating actuators, as considered in the following references \cite{paper1,paper2,paper3,paper4,paper5,paper6,paper7,paper8,paper9,paper10,paper11,paper12,paper13,paper14,paper15,paper16,paper17,paper18,paper19,paper20}.


\bibliographystyle{IEEEtran} 
\bibliography{references}


\end{document}